\documentclass{gtart_h}

\def\ifplaintex{\expandafter\ifx\csname documentclass\endcsname\relax}

\ifplaintex 
\else
\fi

\expandafter\ifx\csname beginpicture\endcsname\relax
\expandafter\ifx\csname documentclass\endcsname\relax
\input pictex \else
\input prepictex \input pictex \input postpictex \fi\fi

\def\gt{{\mathsurround=0pt\it $\cal G\mskip-2mu$eometry \&\ 
$\cal T\!\!$opology}}        

\def\gtp{{\mathsurround=0pt\it $\cal G\mskip-2mu$eometry \&\ 
$\cal T\!\!$opology $\cal P\!$ublications}}  

\def\lognumber#1{\def\thelognumber{#1}}
\def\volumenumber#1{\def\thevolumenumber{#1}}
\def\papernumber#1{\def\thepapernumber{#1}}
\def\volumeyear#1{\def\thevolumeyear{#1}}

\def\pagenumbers#1#2{\def\startpage{#1}\def\finishpage{#2}}
\def\published#1{\def\publishdate{#1}}
\def\proposed#1{\def\theproposer{#1}}
\def\seconded#1{\def\theseconders{#1}}
\def\received#1{\def\receiveddate{#1}}
\def\revised#1{\def\reviseddate{#1}}
\def\accepted#1{\def\accepteddate{#1}}

\def\coverauthors#1{\def\thecoverauthors{#1}}
\def\asciiauthors#1{\def\theasciiauthors{#1}}
\def\asciiaddress#1{\def\theasciiaddress{#1}}
\def\asciiemail#1{\def\theasciiemail{#1}}

\long\def\asciiabstract#1{\long\def\theasciiabstract{#1}}

\let\\\par\let\thelognumber\relax
\let\thevolumenumber\relax\let\thepapernumber\relax
\let\thevolumeyear\relax\let\thesamplenumber\relax\let\startpage\relax
\let\finishpage\relax\let\publishdate\relax\let\receiveddate\relax
\let\reviseddate\relax\let\accepteddate\relax\let\theasciititle\relax
\let\theasciiauthors\relax\let\theasciiaddress\relax
\let\theasciiabstract\relax
\let\theasciiemail\relax\let\theshortauthors\relax\let\theshorttitle\relax
\let\thecoverauthors\relax

\long\def\maketitlep{   

\count0=\startpage

\gt\hfill      
\beginpicture
\setcoordinatesystem units <0.33truein, 0.33truein> point at 2.2 0.9
\setplotsymbol ({$\cal G$})
\plotsymbolspacing=9truept
\circulararc 315 degrees from 0 1 center at 0 0
\setplotsymbol ({$\cal T$})
\circulararc 315 degrees from 1 -1 center at 1 0
\endpicture
\break
{\small\ifx\thesamplenumber\relax 
Volume \else Sample
\fi\thevolumenumber\ (\thevolumeyear)
\startpage--\finishpage\nl
Published: \publishdate}
\vglue 0.5truein plus 0.4fil minus 0.1truein

{\parskip=0pt\leftskip 0pt plus 1fil\def\\{\par\smallskip}{\ifplaintex\large
\else\Large\fi\bf\thetitle}\par\medskip}   

\vglue 0pt plus 0.1fil 

{\parskip=0pt\leftskip 0pt plus 1fil\def\\{\par}{\sc\theauthors}
\par\medskip}

\vglue 0pt plus 0.1fil 

{\small\parskip=0pt\let\newline\\
{\leftskip 0pt plus 1fil\def\\{\par}{\sl\theaddress}\par}
\expandafter\ifx\theemail\relax    
\relax\else\vglue 5pt plus 0.02fil minus 2pt\def\\{\stdspace{\rm 
and}\stdspace} 
\cl{Email:\stdspace\tt\theemail}\fi
\ifx\theurl\relax                  
\relax\else\vglue 5pt plus 0.02fil minus 2pt\def\\{\stdspace{\rm 
and}\stdspace}
\cl{URL:\stdspace\tt\theurl}\fi\par}

\vglue 7pt plus 0.3fil minus 3pt

{\bf Abstract}
\vglue 5pt plus 0.1fil minus 2pt

\theabstract

\vglue 7pt plus 0.3fil minus 3pt

{\bf AMS Classification numbers}\quad Primary:\quad \theprimaryclass

Secondary:\quad \thesecondaryclass

\vglue 5pt plus 0.3fil minus 2pt

{\bf Keywords:}\quad \thekeywords

\vglue 10pt plus 0.5fil minus 5pt

{\small  Proposed: \theproposer\hfill Received: \receiveddate\nl
Seconded: \theseconders\hfill 
\ifx\reviseddate\relax                         
Accepted: \accepteddate                        
\else
Revised: \reviseddate                          
\fi}
\eject
}       

\let\maketitlepage\maketitlep
\let\maketitle\maketitlepage

\font\phead=cmsl9 scaled 950
\font=cmsl9 scaled 1050
\font\pnum=cmbx10 scaled 913
\font\lnum=cmbx10 
\font\pfoot=cmsl9 scaled 950
\font=cmsl9 scaled 1050
\ifplaintex
\headline{\vbox to 0pt{\vskip -4.5mm\line{\small\phead\ifnum
\count0=\startpage ISSN 1364-0380 (on line)
1465-3060 (printed) \hfill {\pnum\folio}\else\ifodd\count0\def\\{ }%
\ifx\theshorttitle\relax\thetitle\else\theshorttitle\fi\hfill{\pnum\folio}
\else\def\\{ and }{\pnum\folio}\hfill\ifx\theshortauthors\relax\theauthors
\else\theshortauthors\fi\fi\fi}\vss}}
\footline{\vbox to 0pt{\vglue 0mm\line{\small\pfoot\ifnum\count0=\startpage
\copyright\ \gtp\hfill\else
\gt, Volume \thevolumenumber\ (\thevolumeyear)\hfill\fi}\vss
}}
\else
\makeatletter
\def\@oddhead{{\small\lhead\ifnum\count0=\startpage ISSN 1364-0380 (on line)
1465-3060 (printed) \hfill {\lnum\number\count0}\else\ifodd\count0
\def\\{ }\ifx\theshorttitle\relax \thetitle \else\theshorttitle\fi\hfill
{\lnum\number\count0}\else\def\\{ and }{\lnum\number\count0}
\hfill\ifx\theshortauthors\relax 
\theauthors\else\theshortauthors\fi\fi\fi}}\def\@evenhead{\@oddhead}
\def\@oddfoot{\small\lfoot\ifnum\count0=\startpage\copyright\ \gtp\hfill\else
\gt, Volume \thevolumenumber\ (\thevolumeyear)\hfill\fi}
\def\@evenfoot{\@oddfoot}
\makeatother
\fi

\newwrite\gtoutfile
\long\gdef\makeheadfile{  
{\def\\{, }\def\s{ }
\immediate\openout\gtoutfile head.xxx
\immediate\write\gtoutfile{Proxy-for: \ifx\theasciiauthors\relax
\theauthors\else\theasciiauthors\fi\s<\ifx\theasciiemail\relax\theemail\else\theasciiemail\fi>}
\immediate\write\gtoutfile{\noexpand\\}
\immediate\write\gtoutfile{Authors: \ifx\theasciiauthors\relax
\theauthors\else\theasciiauthors\fi}
\immediate\write\gtoutfile{Title: \ifx\theasciititle\relax
\thetitle\else\theasciititle\fi}
\immediate\write\gtoutfile{Subj-class: GT or SG or MG etc}
\immediate\write\gtoutfile{MSC-class: \theprimaryclass\ifx\thesecondaryclass\relax\else, \thesecondaryclass\fi}
\immediate\write\gtoutfile{Journal-ref: Geom. Topol. \thevolumenumber
(\thevolumeyear) \startpage-\finishpage}
\immediate\write\gtoutfile{Comments: Published by Geometry and Topology at}
\immediate\write\gtoutfile{\s\s http://www.maths.warwick.ac.uk/gt/GTVol\thevolumenumber/paper\thepapernumber.abs.html}
\immediate\write\gtoutfile{\noexpand\\}
\immediate\write\gtoutfile{}
\ifx\theasciiabstract\relax
\immediate\write\gtoutfile{\theabstract}\else
\immediate\write\gtoutfile{\theasciiabstract}\fi
\immediate\write\gtoutfile{}
\immediate\write\gtoutfile{\noexpand\\}
\immediate\write\gtoutfile{}
\immediate\closeout\gtoutfile}}  

\def\maketitlepage{\maketitlep\makeheadfile}
\let\maketitle\maketitlepage

\lognumber{392}
\volumenumber{8}\papernumber{7}\volumeyear{2004}
\pagenumbers{295}{310}
\received{7 December 2003}
\revised{9 December 2003}
\published{14 February 2004}
\accepted{13 February 2004}
\proposed{Robion Kirby}
\seconded{John Morgan, Ronald Stern}

\usepackage{amsmath,amssymb}

\newcommand{\rem}{Remark}
\newcommand{\rems}{Remarks}
\newcommand{\ack}{Acknowledgements}

\theoremstyle{plain}
\newtheorem{thm}{Theorem}
\newtheorem{cor}[thm]{Corollary}
\newtheorem{conj}[thm]{Conjecture}
\newtheorem{lem}[thm]{Lemma}
\newtheorem{prop}[thm]{Proposition}
\theoremstyle{definition}

\newcommand{\spinc}{\mathrm{Spin}^c}
\newcommand{\s}{\mathfrak{s}}
\newcommand{\R}{\mathbb{R}}

\newcommand{\Q}{\mathbb{Q}}
\newcommand{\CP}{\mathbb{CP}}
\newcommand{\Z}{\mathbb{Z}}

\newcommand{\SO}{\mathit{SO}}

\newcommand{\A}{\mathbb{A}}
\newcommand{\HF}{\mathit{HF}}
\newcommand{\HFF}{\mathit{HFF}}
\makeatletter
\@ifpackageloaded{mathpi}{%
   \newcommand{\DS}{{\mathcal{D}\!}}
   }{\newcommand{\DS}{\mathcal{D}}}
\makeatother
\newcommand{\SW}{\mathit{SW}}
\newcommand{\bD}{\bar{D}}
\newcommand{\Hol}{\mathrm{Hol}}

\title          {Witten's conjecture and Property P}

\authors        {P\,B Kronheimer\\T\,S Mrowka}
\coverauthors   {P\noexpand\thinspace B Kronheimer\\T\noexpand\thinspace S Mrowka}
\asciiauthors   {P B Kronheimer and T S Mrowka}
\address        {Department of Mathematics, Harvard University\\Cambridge 
                 MA 02138, USA}
\secondaddress  {Department of Mathematics, Massachusetts Institute of 
                 Technology\\Cambridge MA 02139, USA} 
\gtemail         {\mailto{kronheim@math.harvard.edu}\qua{\rm and}\qua
                  \mailto{mrowka@math.mit.edu}}
\asciiaddress   {Department of Mathematics, Harvard University\\Cambridge 
                 MA 02138, USA\\and\\Department of Mathematics, Massachusetts 
                 Institute of Technology, Cambridge MA 02139, USA}
\asciiemail     {kronheim@math.harvard.edu, mrowka@math.mit.edu}


\begin{abstract}
       Let $K$ be a non-trivial knot in the $3$--sphere and let $Y$ be the
        $3$--manifold obtained by surgery on $K$ with surgery-coefficient
        $1$. Using tools from gauge theory and symplectic topology, it
        is shown that the fundamental group of $Y$ admits a non-trivial
        homomorphism to the group $\SO(3)$. In particular, $Y$ cannot be a
        homotopy-sphere.
\end{abstract}

\asciiabstract{%
        Let K be a non-trivial knot in the 3-sphere and let Y be the
        3-manifold obtained by surgery on K with surgery-coefficient
        1. Using tools from gauge theory and symplectic topology, it
        is shown that the fundamental group of Y admits a non-trivial
        homomorphism to the group SO(3). In particular, Y cannot be a
        homotopy-sphere.}

\begin{document}
\maketitlepage

\section{Introduction}

Let $K$ be a knot in $S^{3}$ and let $Y_{1}$ be the oriented 3--manifold
obtained by $+1$--surgery on $K$. The following is one
formulation of the ``Property P'' conjecture for knots:

\begin{conj}\label{conj:PP}
    If $K$ is a non-trivial knot, then $Y_{1}$ is not a homotopy
    $3$--sphere.
\end{conj}

The purpose of this note is to prove the conjecture. The ingredients
of the argument are: (a) Taubes' theorem \cite{Taubes1} on the
non-vanishing of the Seiberg--Witten invariants for symplectic
$4$--manifolds; (b) the theorem of Gabai \cite{Gabai3} on the existence
of taut foliations on $3$--manifolds with non-zero Betti number; (c)
the construction of Eliashberg and Thurston
\cite{Eliashberg-Thurston}, which produces a contact structure from a
foliation; (d) Floer's exact triangle \cite{Floer2, Braam-Donaldson}
for instanton Floer homology; (e) a recent result of Eliashberg
\cite{Eliashberg} on concave
filling of contact $3$--manifolds\footnote{%
The authors have learned that this result
was also known to Etnyre, who shows in \cite{Etnyre} that it is a
straightforward extension of the earlier results of
\cite{Etnyre-Honda}.}; and (f) Witten's conjecture relating
the Seiberg--Witten and Donaldson invariants of smooth $4$--manifolds.
Although the full version of Witten's conjecture remains open, a
weaker version that is still strong enough to serve our purposes has
recently been established by Feehan and Leness \cite{Feehan-Leness},
following a program proposed by Pidstrigatch and Tyurin.  With these
ingredients, we shall prove:

\begin{thm}\label{thm:main}
    Let $Y_{1}$ be obtained by $+1$--surgery on a non-trivial knot $K$
    in $S^{3}$. Then there is a non-trivial homomorphism
    $\rho\co\pi_{1}(Y_{1}) \to \SO(3)$.
\end{thm}

It is known \cite{Gordon-Luecke} that surgery on a non-trivial knot can never yield
$S^{3}$, so the Property P conjecture would
follow from the Poincar\'e conjecture. Theorem~\ref{thm:main} is a
slightly sharper statement which implies Conjecture~\ref{conj:PP}. The
same techniques yield a closely-related theorem:

\begin{thm}\label{thm:Betti}
    Let $Y$ be an irreducible, closed, orientable
    $3$--manifold (not $S^{1}\times S^{2}$), and
    let $v$ be an element of
    $H^{2}(Y;\Z/2)$. Then 
    there is a homomorphism $\rho \co 
    \pi_{1}(Y) \to \SO(3)$  having $w_{2}(\rho)
    = v$.
\end{thm}

\subparagraph{\rems}
The question whether surgery on a knot could produce a counterexample to
the Poincar\'e conjecture was asked explicitly by Bing in \cite{Bing-Erratum},
and the question was formalized with the definition of
``Property P'' by Bing and Martin in
\cite{Bing-Martin}. To verify that a knot $K$ has Property P in their
sense, it is sufficient to verify that the $3$--manifolds $Y$ obtained by
non-trivial Dehn surgeries on $K$ all have non-trivial fundamental group.
In \cite{CGLS}, it was shown
that $\pi_{1}(Y)$ is non-trivial if $K$ is non-trivial and the
surgery-coefficient is not $\pm 1$. This is why
Conjecture~\ref{conj:PP} is now equivalent to the original version.
The problem appears on Kirby's problem list \cite[Problem
1.15]{Kirby}, where there is also a summary of some of the contributions
that have been made.

It follows from Casson's work (see \cite{Akbulut-McCarthy}) that if
$Y$ is obtained by Dehn surgery on a knot $K$ whose symmetrized
Alexander polynomial $\Delta_{K}$ satisifes $\Delta''_{K}(1) \ne 0$,
then $\pi_{1}(Y)$ admits a non-trivial homomorphisms to $\SO(3)$. Such
knots therefore have Property P.  
The argument used by Casson is closely related to what is
done here. The quantity $\Delta''_{K}(1)$ is equal to the Euler
characteristic of a Floer homology group $\HF(Y_{0})$, associated to
the manifold $Y_{0}$
obtained by $0$--framed surgery on $K$. We shall show that the Floer
homology group $\HF(Y_{0})$ itself is always non-trivial if $K$ is not
the unknot, even though the Euler characteristic may vanish. 

The authors were aware some time ago that Property P could be deduced
from Witten's conjecture and other known results,
if one only had a suitably general ``concave
filling'' result for symplectic $4$--manifolds with contact boundary,
as explained later in this paper.
At the time (around 1996), no concave filling results were known.
The first general result on concave filling of contact $3$--manifolds
is given in \cite{Etnyre-Honda}, using results on open-book
decompositions from \cite{Giroux}. More recently, Eliashberg has shown
\cite{Eliashberg} that one can construct a concave filling compatible
with a given symplectic form on a collar of the contact $3$--manifold,
provided only that the symplectic form is positive on the contact
planes. It is this stronger result from \cite{Eliashberg} that we
need here.

\subparagraph{\ack}
The first author was supported by NSF grant
DMS-0100771. The second author was supported by NSF grants 
DMS-0206485, DMS-0111298 and FRG-0244663. Both authors would like to
thank Yasha Eliashberg for generously sharing his expertise.

\section{Donaldson and Seiberg--Witten invariants}

\subsection{Donaldson invariants and simple type}

Let $X$ be a smooth, closed, oriented $4$--manifold, with $b^{+}(X)$
odd and greater than $1$, and $b_{1}(X) = 0$. Fix a homology
orientation for $X$. For each $w \in
H^{2}(X;\Z)$, the Donaldson invariants of $X$ (constructed using
$U(2)$ bundles with first Chern class $w$) constitute a linear map
\[
       D_{X}^{w} \co  \A(X) \to  \Z,
\]
where $\A(X)$ is the symmetric algebra on $H_{2}(X;\Z) \oplus
H_{0}(X;\Z)$. Our notation here follows \cite{KM-Structure}, and we
write $x$ for the element of $\A(X)$ corresponding to the positive
generator of $H_{0}(X;\Z)$. We make $\A(X)$ a graded algebra, by
putting the generators from $H_{2}(X;\Z)$ in degree $2$ and the
generator $x$ in degree $4$. With this grading, the restriction
\[
        D_{X}^{w} \co  \A_{2d}(X) \to  \Z
\]
is non-zero only when
\begin{equation}\label{eq:d-congruence}
        d \equiv -w^{2}-\frac{3}{2} (b^{+}(X) + 1) \pmod {4}.
\end{equation}

The manifold $X$ is said to have
\emph{simple type} if the invariant satisfies
\[
    D_{X}^{w}(x^{2} z) = 4 D_{X}^{w}(z)
\]
for all $z$ in $\A(X)$.  This notion was introduced in
\cite{KM-Structure}, where it was shown that $X$ has simple type if it
contains a \emph{tight surface}: a smoothly embedded oriented surface $\Sigma$
whose genus $g$ satisfies $2g-2 = [\Sigma]\cdot[\Sigma] > 0$. For
manifolds of simple type, it is natural to introduce
\[
           \bD^{w}_{X} \co  \A(X) \to \Z,
\]
$$ 
\bD^{w}_{X}(z) = D^{w}_{X}(z)  + D^{w}_{X}(zx/2). \leqno{\hbox{defined by}}
$$
If $z$ is homogeneous of degree $2d$, then only one of the terms on
the right can be non-zero because of the congruence
\eqref{eq:d-congruence}; and both terms are zero unless
\begin{equation}\label{eq:d-congruence2}
        d \equiv -w^{2}-\frac{3}{2} (b^{+}(X) + 1) \pmod {2}.
\end{equation}
We combine the Donaldson invariants to form
a series
 \[
   \begin{aligned}
   \DS_{X}^w(h)\ &=\  \bD^w_X(e^h) \\
	    &=\ \sum  D^w_X(h^d)/d!\ +\ \tfrac{1}{2} \sum  D^w_X(x
            h^d)/d!
   \end{aligned}            
 \]
We regard this as a formal power series for $h\in H_{2}(X;\R)$. The
main result of \cite{KM-Structure} contains the following:

\begin{thm}\label{thm:struct}
{\rm\cite{KM-Structure}}\qua
Let $X$ be a 4--manifold of simple type with $b_{1}=0$.
Then the
Donaldson series converges for all $w$ and
there exist finitely many cohomology classes $K_1$, \dots, $K_s\in
H^2(X;\Z)$ and non-zero rational numbers $\beta_1$, \dots, $\beta_s$
(both independent of $w$) such that
 \[
	\DS_{X}^w\ =\ \exp\left(\frac{Q}{2}\right) \,\sum_{r=1}^s 
		(-1)^{(w^2 + K_r\cdot w)/2}\,\beta_r e^{K_r}
 \]
as analytic functions on $H_2(X;\R)$.  Here $Q$ is the intersection
form, regarded as a quadratic function.  Each of the classes $K_r$
is an integral lift of $w_2(X)$.  
\end{thm}

\subparagraph{\rems} The classes $K_{r}$ are called the \emph{basic
classes} of $X$. The theorem is supposed to include the case
that the Donaldson invariants are identically zero. This is the case
$s=0$. The Donaldson series is always either an even or an odd
function of $h$, so the non-zero basic classes come in pairs
differing by sign.

It is more common today to use the terms ``simple type'' and ``basic
classes'' to refer to properties defined not by the Donaldson
invariants but by the Seiberg--Witten invariants, as explained below.
We will therefore refer to these as \emph{$D$--simple type} and
\emph{$D$--basic classes} henceforth, to avoid ambiguity.

\subsection{Seiberg--Witten invariants and Witten's conjecture}

The Seiberg--Witten invariants of a $4$--manifold $X$ such as the one we
are considering (with $b^{+}$ odd and greater than $1$ and $b_{1}=0$)
are a function on the set of $\spinc$ structures on $X$. For each $\spinc$
structure $\s$, they define an integer
$
        \SW(\s) \in \Z.
$
To simplify our notation, we shall assume that $X$ has no $2$--torsion
in its second cohomology: in this case, $\s$ is determined by the
first Chern class $K$ of the corresponding half-spin bundle $S^{+}$,
and we can regard $\SW$ as a function of $K$:
\[
        \SW \co  H^{2}(X;\Z) \to \Z.
\]
The manifold $X$ is said to have \emph{$\SW$--simple type} if
$\SW(K)=0$
whenever $K^{2}$ is not equal to $2\chi + 3\sigma$.  The
\emph{$\SW$--basic classes} are the classes $K \in H^{2}(X;\Z)$ with
$\SW(K)\ne 0$. The following is a stripped-down version of Witten's
conjecture from \cite{Witten}.

\begin{conj}\label{conj:Witten}
    Let $X$ be a $4$--manifold with $b^{+}$ odd and greater than $1$,
    with $b^{1}(X)=0$ and with no $2$--torsion in $H^{2}(X;\Z)$.
    Suppose $X$ has $\SW$--simple type. Then $X$ has $D$--simple type,
    the $D$--basic classes are the $\SW$--basic classes, and for each
    basic class $K_{r}$, the corresponding rational number $\beta_{r}$
    in the statement of Theorem~\ref{thm:struct} is given by
    \[
                \beta_{r} = c(X)
                \SW(K_{r}),
    \]
    where $c(X)$ is a non-zero rational number depending on $X$.
\end{conj}

An important corollary of this conjecture is the assertion that the Donaldson
invariants are non-zero if the Seiberg--Witten invariants are non-zero
and $X$ has $\SW$--simple type. Witten's conjecture also gives the
value of $c(X)$ as
\[
      c(X) = 2^{2 + \frac{1}{4}(7\chi + 11\sigma)},
\]
but we will not need this statement.

\subsection{The theorem of Feehan and Leness}

A weaker version of Witten's conjecture is proved by Feehan and
Leness in
\cite{Feehan-Leness}. We rephrase Theorem~1.1 of \cite{Feehan-Leness}
here, specializing to the case that $X$ has $\SW$--simple type, and
simplifying the statement to suit our needs, as follows.  The theorem involves a
choice of auxiliary class $\Lambda \in H^{2}(X;\Z)$ with $\Lambda - w
= w_{2}(X)$ mod $2$.  In the version we state here, we take $\Lambda$
to be the class dual to a tight surface in $X$. This ensures that
$\Lambda\cdot K$ is zero, for all $\SW$--basic classes $K$. The
presence of a tight surface
ensures that $X$ has $D$--simple type. We choose $\Lambda$ to be
divisible by $2$ and $w$ to be an integer lift of $w_{2}(X)$. Set
\[
          N = \Lambda^{2} \in \Z.
\]
We may replace $\Lambda$ by any multiple of $\Lambda$, to make $N$ as
large as we might need.

\begin{thm}\label{thm:FL}{\rm\cite{Feehan-Leness}}\qua
    Let $X$ be a $4$--manifold with $b_{1}=0$ and $b^{+}$ odd and
    greater than $1$. Suppose that $X$ has no $2$--torsion in its
    second cohomology and has $\SW$--simple type. Suppose in addition
    that $X$ contains a tight surface with positive self-intersection
    number. Let $\Lambda$ and $N$
    be as above, and let $d$ be an integer in the range
    \[
            0 \le d < N  - \tfrac{1}{4}(\chi + \sigma) - 2
    \]
    satisfying the congruence \eqref{eq:d-congruence2}.
    Then for any class $h$ in $H_{2}(X;\R)$ with $\Lambda\cdot h  =
    0$, we have
    \[
           \bD^{w}_{X}(h^{d}) = \sum_{K}
             (-1)^{(w^2 + K\cdot w)/2} \SW(K)
            \,
                    p_{d}(K\cdot h, Q(h)).
    \]
    Here $p_{d}$ is a weighted-homogeneous polynomial,
    \[
                  p_{d}(s,t) = \sum_{a + 2b = d} C_{a,b} s^{a} t^{b},
    \]
    whose coefficients $C_{a,b}\in \Q$ are universal functions of $\chi(X)$,
    $\sigma(X)$ and $N$.
\end{thm}

From this result, it is straightforward to deduce:

\begin{cor}\label{cor:tightWitten}
    Witten's conjecture, in the form of Conjecture~\ref{conj:Witten},
    holds for $X$ as long as $X$ satisfies the following three
    additional conditions:
    \begin{enumerate}
        \item $X$ contains a tight surface $\Sigma$ with positive
        self-intersection;
        \item $X$ has the same Euler number and signature as some smooth hypersurface
        in $\CP^{3}$ whose degree is even and at least 6;
        \item $X$ contains a sphere of self-intersection $-1$.
    \end{enumerate}
\end{cor}

\subparagraph{\rem} The second condition is much more restrictive
than necessary, but suffices for our application.

\begin{proof}[Proof of the corollary]
    The assertion of Conjecture~\ref{conj:Witten} is an equality
    \begin{equation}\label{eq:witten}
               \DS^{w}_{X} = c(X)\exp(Q/2) \sum_{K} (-1)^{(w^2 + K\cdot w)/2}
                         \SW(K) e^{K}    
    \end{equation}
    of analytic functions on $H^{2}(X;\R)$, where $c(X)$ is a non-zero
    rational number.  We are assuming that $X$ contains a tight
    surface, so $X$ has $D$--simple type.
    If we change $w$ to $w'$, then we know how
    $\DS_{X}^{w}$ changes, from Theorem~\ref{thm:struct}, and we know also
    how the right-hand side changes. It is therefore enough to verify
    the conjecture for one particular $w$. We take $w$ to be an
    integral lift of $w_{2}(X)$.


    Let $\Lambda \in H^{2}(X;\Z)$ be some large even multiple of the
    class dual to the tight surface.  All the $\SW$--basic classes
    and all the $D$--basic classes are orthogonal to $\Lambda$ by the
    adjunction inequality. If we write $h= h_{1} + h_{2}$, where
    $\Lambda\cdot h_{1} = 0$ and $h_{2}$ is in the span of the dual of
    $\Lambda$, then
    \[
              \DS_{X}^{w}(h_{1} + h_{2}) = \DS_{X}^{w}(h_{1}) \exp(
              Q(h_{2})/2).
    \]
    The same holds for the function defined by the right-hand side of
    \eqref{eq:witten}. So it is enough to verify that the conjecture
    holds for the restriction of the Donaldson series to the kernel of
    $\Lambda$.

    Let $X_{*}$ be a hypersurface in $\CP^{3}$ with  the same Euler number
    and signature as $X$. We take $\Lambda_{*}$ in $H^{2}(X_{*};\Z)$ to be a
    class orthogonal to the canonical
    class $K_{*}$ of $X_{*}$, represented by a tight surface.
    By replacing $\Lambda_{*}$ and $\Lambda$ by
    suitable multiples, we can arrange that they have the same square
    $N$.  When the degree of $X_{*}$ is even, the congruence
    \eqref{eq:d-congruence2} asserts that $d$ is even. The Donaldson
    invariants of $X$ and $X_{*}$ are even functions on the second
    homology in this case, and the Seiberg--Witten invariants satisfy
    $\SW(K)= \SW(-K)$ in both cases.

    We apply Theorem~\ref{thm:FL} to $X_{*}$, with $w=0$. The $\SW$--basic classes
    are $\pm K_{*}$, and  it is known that $\SW(\pm K_{*}) = 1$.
    We learn that
    \[
             \bD^{0}_{X_{*}}(h^{d}) = 2 \sum_{a+2b = d}
             C_{a,b}(K_{*}\cdot h)^{a} Q(h)^{b}
    \]
   for all $h$ orthogonal to $\Lambda_{*}$. This formula determines
   the coefficients $C_{a,b}$ entirely, in terms of the Donaldson
   invariants of $X_{*}$, because the linear function $K_{*}$ and the
   quadratic form $Q$ are algebraically independent as functions on
   this vector space.

   In particular, we see that $C_{a,b}$ is
   independent of $N$.
   We can therefore sum over all $d$, and write
   \[
     \DS^{0}_{X_{*}}(h) = 2 \sum_{\text{$d$ even}}(1/d!) \sum_{a+2b=d}
     C_{a,b}(K_{*}\cdot h)^{a} Q(h)^{b}.
   \]
   On the other hand, we know from Theorem~\ref{thm:struct} that
   $\DS^{0}_{X_{*}}$ has the special form given there; and we also
   know that the Donaldson invariants of this complex surface are not
   identically zero. Thus
   \[
        2 \sum_{\text{$d$ even}}(1/d!) \sum_{a+2b=d}
     C_{a,b}(K_{*}\cdot h)^{a} Q(h)^{b} = \exp(Q(h)/2) f(K_{*}\cdot h)
   \]
   where $f\co \R\to\R$ is a non-zero even function of the form
   \[
             f(t) = \sum_{r=1}^{m} \alpha_{r} \cosh (\lambda_{r} t)
   \]
   for some rational numbers $\alpha_{r}$ and $\lambda_{r}\ge 0$. The
   rational numbers $\lambda_{r}$ are such that the $D$--basic classes
   of $X$ are $\pm \lambda_{r} K_{*}$. The basic classes are supposed
   to be integer classes, and this constrains the denominator of
   $\lambda_{r}$. The adjunction inequality also implies that
   $\lambda_{r}\le 1$.

   With this information about $C_{a,b}$, we can now apply
   Theorem~\ref{thm:FL} to our original $X$, to deduce that
   \[
       \DS^{w}_{X} = \exp(Q/2) \sum_{K} (-1)^{(w^2 + K\cdot w)/2}
                         \SW(K) f(K)    
   \]
   as functions on the orthogonal complement of $\Lambda$. If any of
   the $\lambda_{r}$ are not integral, then this formula is
   inconsistent with Theorem~\ref{thm:struct}, because the $\SW$--basic
   classes $K$ for $X$ are primitive, because $X$ contains a sphere of
   square $-1$. The $D$--basic classes are also all non-zero for $X$,
   for the same reason, and this means that no $\lambda_{r}$ can be
   zero. So $\lambda_{r}$ can only be $\pm 1$, and it follows that
   $f(K)$ is simply a multiple of $\cosh(K)$. This establishes the
   result.
\end{proof}

\section{Proofs of the theorems}
\label{sec:Outline}

\subsection{Concave filling}

Let $Y$ be a closed oriented $3$--manifold (not necessarily connected),
and $\xi$ an oriented contact
structure compatible with the orientation of $Y$. This means that
$\xi$ is the $2$--plane field defined as the kernel of a $1$--form
$\alpha$ on $Y$, and $\alpha\wedge d\alpha$ is a positive $3$--form. If
$Y$ is the oriented boundary of an oriented $4$--manifold of $W$, then
a symplectic form $\omega$ on $W$ is said to be weakly compatible
with $\xi$ if the restriction $\omega|_{Y}$ is positive on the
$2$--plane field $\xi$; or equivalently, if $\alpha\wedge \omega|_{Y} >
0$.  
The following is proved in \cite{Eliashberg}.

\begin{thm}\label{thm:Eliash}{\rm\cite{Eliashberg}}\qua
Let $Y$ be the oriented boundary of a $4$--manifold $W$ and let $\omega$ be a
symplectic form on $W$. Suppose there is a contact structure $\xi$ on
$Y$ that is compatible with the orientation of $Y$ and weakly
compatible with $\omega$. Then we can embed $W$ in a closed symplectic
$4$--manifold $(X,\Omega)$ in such a way that $\Omega|_{W} = \omega$.
\end{thm}

Because we will need to construct an $(X,\Omega)$ satisfying some
additional mild restrictions, we summarize how $X$ is constructed
in \cite{Eliashberg} as
a smooth manifold (without concern for the symplectic form). If the
components of $Y$ are $Y_{1},\dots,Y_{n}$, then the first step is to
choose an open-book decomposition of each $Y_{i}$ with binding
$B_{i}$. These open-book decompositions are required to be compatible
with the contact structures $\xi|_{Y_{i}}$ in the sense of
\cite{Giroux}. We can take each binding $B_{i}$ to be connected. Let
$W'$ be obtained from $W$ by attaching a $2$--handle along each knot
$B_{i}$ with zero framing. The boundary $Y' = \partial W'$ is the union of
$3$--manifolds $Y_{i}'$, obtained from $Y_{i}$ by zero surgery: each
$Y_{i}'$ fibers over the circle with typical fiber $\Sigma_{i}$. The
genus of $\Sigma_{i}$ is the genus of the leaves of the open-book
decomposition of $Y_{i}$. For
each $i$, one then constructs a symplectic Lefschetz fibration
\begin{equation}\label{eq:Lef-i}
          p_{i} \co  Z_{i} \to B_{i}
\end{equation}
over a $2$--manifold-with-boundary $B_{i}$, with $\partial B_{i}
=S^{1}$. One constructs $Z_{i}$ to have
the same fiber $\Sigma_{i}$, and $\partial Z_{i} = - Y'_{i}$. The
$4$--manifold $X$ is obtained as the union of $W'$ and the $Z_{i}$,
joined along their common boundaries $Y_{i}'$.

There is considerable freedom in this construction. We exploit this
freedom in a sequence of lemmas, each of which states that we can
choose $Z_{i}$ so as to fulfill a particular additional property.

\begin{lem}
    \label{lemma:Flux}
    We can choose the Lefschetz fibration $p_{i} \co  Z_{i} \to B_{i}$ so
    that the base $B_{i}$ is a disk $D^{2}$.
\end{lem}

\begin{proof}
    The constructions in \cite{Eliashberg} already establish this. We
    present a slight variation of the argument.

    As a component of
    $\partial W'$, the $3$--manifold $Y_{i}'$ carries a $2$--form
    $\eta' \in \Omega^{2}(Y_{i}')$,
    the restriction of the symplectic form $\omega'$ from $W'$. This
    form is positive on the fibers of the fibration $p'\co  Y_{i}' \to
    S^{1}$, and its kernel is a line-field on $Y_{i}'$ transverse to
    the fibers. There is a unique vector field $V'$ on $Y_{i}'$
    contained in the line-field, with $p'_{*}(V') = \partial/\partial
    S^{1}$ on the circle $S^{1}$. The flow generated by $V'$ preserves
    $\eta'$; and at time $2\pi$ the flow determines a holonomy
    automorphism $\Hol(\eta')$ of the fiber over $1\in S^{1}$, which
    is an area-preserving map of the surface.

    Since positive Dehn twists generate the mapping class group, we
    can construct a Lefschetz fibration $p^{0}_{i} \co  Z^{0}_{i} \to D^{2}$
    whose boundary is topologically $-Y_{i}'$, as a surface bundle
    over $S^{1}$. This Lefschetz fibration can be made symplectic; and
    we write $\eta''$ for the restriction of the symplectic form from
    $Z^{0}_{i}$ to
    $Y_{i}'$. We can assume that $\eta'$ and $\eta''$ have the same
    integral on the fiber $\Sigma_{i}$.

    If we can choose $Z^{0}_{i}$ so that
    \begin{equation}\label{eq:Hol-Hol}
               \Hol(\eta') = \Hol(\eta'')
    \end{equation}
    as area-preserving maps of the fiber over $1$, then there is a
    fiber-preserving diffeomorphism $\psi$ of $Y_{i}'$ with
    $\psi^{-1}(\eta')=\eta''$. We can then use $\psi$ to
    attach $Z^{0}_{i}$ to $W'$
    along $Y_{i}'$ (see \cite{Eliashberg}) and our task will be
    complete.
    At this point however, we only know that the map
    $\phi=\Hol(\eta') \circ
    \Hol(\eta'')^{-1}$ is isotopic to the identity in
    $\mathrm{Diff}(\Sigma_{i})$.  

    To complete the proof of
    the lemma, it will be enough to construct a symplectic Lefschetz fibration
    \[
            p \co  (V,\omega_{V}) \to D^{2}
    \]
    whose boundary is the topologically trivial surface bundle over
    $S^{1}$ and
    whose holonomy is given by $\Hol(\eta) = \phi$, where $\eta =
    \omega_{V}|_{\partial V}$. We can then form
    $Z_{i}$ as the union of $Z_{i}^{0}$ and $V$, attached along a
    neighborhood of a fiber in their boundaries. That such a $V$ exists
    is the content of the next lemma, which is a variant of
    \cite[Lemma 3.4]{Eliashberg}.
\end{proof}

\begin{lem}
   Let $\Sigma$ be a closed symplectic surface of area $1$ and genus
   $2$ or more.
    Let $\phi\co  \Sigma \to \Sigma$ be an area-preserving map that is 
    isotopic to the identity through diffeomorphisms. Then there is a
    symplectic Lefschetz fibration $p \co  (V,\omega) \to D^{2}$ with
    $p^{-1}(1) = \Sigma$ and
    $\Hol(\omega|_{V}) = \phi$.
\end{lem}

\begin{proof}
    As explained in \cite{Eliashberg}, it will be enough if we can
    find a $(V,\omega)$ such that $\Hol(\omega_{V})$ has the same
    \emph{flux} as $\phi$.  In this context, the flux has the
    following interpretation. Because the identity component of the
    diffeomorphism group is contractible, we can identify $\partial V$
    with $S^{1}\times \Sigma$ canonically up to fiber-preserving
    isotopy; so we have a canonical map
    \[
               H_{1}(\Sigma) \to H_{2}(\partial V)
    \]
    given by $[\gamma]\mapsto [S^{1}\times \gamma]$.  The flux is the
    element of $H^{1}(\Sigma;\R)$ corresponding to the
    homomorphism
    \[
                   \begin{gathered}
                    f\co  H_{1}(\Sigma) \to \R \\
                    [\gamma] \mapsto \int_{S^{1}\times \gamma}
                    \omega|_{\partial V}.
                   \end{gathered}
    \]
    So the assertion of the lemma is that we can choose
    $p\co (V,\omega) \to D^{2}$ so that the cohomology class of
    $\omega|_{\partial V}$ is any given class in $H^{2}(S^{1}\times
    \Sigma;\R)$, subject only to the constraint that the area of $\Sigma$ is
    $1$.

    To see that this is possible, we observe that we can find first
    an example
    $p_{0} \co  (V_{0},\omega_{0}) \to D^{2}$ whose flux $f$ is
    \emph{zero} and
    such that the map $H_{2}(\partial V_{0};\R) \to H_{2}(V_{0};\R)$
    induced by the inclusion $\partial V_{0} \hookrightarrow V_{0}$ is
    injective.  Such an example is obtained by removing a neighborhood
    of a fiber in a closed Lefschetz fibration $\bar{p}_{0} \co 
    (\bar{V}_{0},\bar{\omega}_{0}) \to S^{2}$; the condition on the
    second homology is achieved if $H_{1}(\bar{V}_{0})$ is zero.

    Next, because non-degeneracy is an open condition on $2$--forms,
    there exists a neighborhood $\mathcal{U}$ of 
    $0\in H^{1}(\Sigma;\R)$ such that, for all $f\in \mathcal{U}$,
    there exists a form $\omega_{f}$ on $V_{0}$ such that
    \[
        p_{0} \co  (V_{0}, \omega_{f}) \to D^{2}
    \]
    is a symplectic Lefschetz fibration whose holonomy on the boundary
    has flux $f$.  Finally, given a general $f$, we can find an
    integer $N$ such that $f/N$ belongs to $\mathcal{U}$. We then
    construct $(V,\omega)$ by attaching $N$ copies of $(V_{0},
    \omega_{f/N})$ along neighborhoods of fibers in their boundaries.
\end{proof}

From now on, we may assume that the base of the fibration $Z_{i}$ is a
disk. We can now
arrange that $H_{1}(Z_{i};\Z)$ is zero. Indeed, $H_{1}(Z_{i};\Z)$ is
generated by a collection of $1$--cycles on the fiber $\Sigma_{i}$, and
we can arrange that these are vanishing cycles in the Lefschetz
fibration.  Thus we can state:

\begin{lem}\label{lem:noH1}
    If the map $H_{1}(Y;\Z) \to H_{1}(W;\Z)$ is surjective, then 
    we can choose $X$ in Theorem~\ref{thm:Eliash} so that
    $H_{1}(X;\Z)$ is zero.
\end{lem}

\begin{proof}
    The hypothesis implies that $H_{1}(Y';\Z) \to H_{1}(W';\Z)$ is
    surjective also. Choose the $Z_{i}$ to have trivial first
    homology, as explained above, and the lemma then follows from the
    Mayer--Vietoris sequence.
\end{proof}

In a similar vein, we have:

\begin{lem}\label{lem:H2onto}
     We can choose $X$ so the restriction map $H^{2}(X;\Z) \to
    H^{2}(W;\Z)$ is surjective.
\end{lem}

\begin{proof}
   The restriction map $H^{2}(W';\Z) \to H^{2}(W;\Z)$ is surjective,
   so we may replace $W$ by $W'$ in the statement.
   If we arrange that $H_{1}(Z_{i};\Z)$ is zero, then the restriction
   map $H^{2}(Z_{i};\Z) \to H^{2}(Y'_{i};\Z)$ is surjective. The
   surjectivity of the map $H^{2}(X;\Z) \to
    H^{2}(W';\Z)$ now follows from the Mayer--Vietoris sequence for
    cohomology.
\end{proof}

We can also specify the Euler number and signature quite freely subject to
some inequalities:

\begin{lem}\label{lem:LikeHyp}
    We can choose $X$ so that its Euler number and signature are the
    same as those of $X_{*}$,
    where $X_{*}$ is a smooth hypersurface in $\CP^{3}$ whose degree
    is even and at least $6$. At the same time, we can arrange that
    $X$ contains a sphere with self-intersection $-1$.
\end{lem}

\begin{proof}
    Our strategy is to arrange that $X$ has the same $b^{+}$ as some
    $X_{*}$ but has smaller $b^{-}$. We then blow up $X$ at enough
    points to make the value of $b^{-}$ agree also.
    
    Let $V\to \CP^{1}$ be a symplectic Lefschetz fibration with
    $b_{1}(V)=0$ and the same fiber genus as $Z_{i}$. Replace $Z_{i}$
    by $\tilde Z_{i}$, the Gompf fiber-sum of $Z_{i}$ and $V$. The
    effect on $b^{+}(X)$ is to add to it the quantity \[ n^{+}(V) = b^{+}(V)
    + 2g-1,\] while $b^{-}(X)$ changes by
    \[
         n^{-}(V) = b^{-}(V)
    + 2g-1.
    \]
    Here $g$ is the fiber genus.  If we use two different
    Lefschetz fibrations, $V$ and $\tilde V$, for which $n^{+}(V)$ and
    $n^{+}(\tilde V)$ are coprime, then the set of values that we can
    achieve for $b^{+}(X)$ includes all sufficiently large integers.

    For hypersurfaces $X_{*}$ in $\CP^{3}$ of large degree, the ratio
    $b^{-}(X_{*}) / b^{+}(X_{*})$ approaches $2$. We can therefore
    achieve our objective by forming a fiber-sum with many copies of
    $V$, provided the ratio $n^{-}(V)/ n^{+}(V)$ satisfies
    \[
        n^{-}(V) / n^{+}(V) < 2.
    \]
    This ratio condition is quite common for Lefschetz fibrations. For
    example, if $S$ is an algebraic surface with an ample class $H$
    satisfying $K_{S} \cdot H > 0$, then the Lefschetz fibration $V$
    constructed from a pencil in the linear system $|dH|$ satisfies
    this inequality, once $d$ is sufficiently large. A $V$ constructed
    in this way may not have the same fiber genus as one of the
    $Z_{i}$, but we can always increase the fiber genus of $Z_{i}$ by
    any positive integer, by
    adjusting the original open-book decomposition of $Y_{i}$.
\end{proof}

We need one last lemma of this sort.

\begin{lem}\label{lem:Xtight}
    We can choose $X$ so that it contains a tight surface of positive
    self-intersection number.
\end{lem}

\begin{proof}
    We can choose a Lefschetz fibration $V\to\CP^{1}$ containing a
    tight surface disjoint from a fiber. We then replace one $Z_{i}$ by
    a Gompf fiber-sum, as in the previous lemma.
\end{proof}

We now combine the conclusions of the last four lemmas with the
construction of Eliashberg and Thurston from
\cite{Eliashberg-Thurston}, to prove the next proposition.

\begin{prop}\label{prop:Summary}
    Let $Y$ be a closed orientable $3$--manifold admitting an oriented taut
    foliation. Suppose $Y$ is not $S^{1}\times S^{2}$. Then $Y$ can be
    embedded as a separating hypersurface in a closed symplectic $4$--manifold
    $(X,\Omega)$. Moreover, we can arrange that $X$ satisfies the
    following additional conditions.
    \begin{enumerate}
        \item The first homology $H_{1}(X;\Z)$ vanishes.

        \item The Euler number and signature of $X$ are the same as
        those of some smooth hypersurface in
        $\CP^{3}$, whose degree is even and not less than $6$.

         \item The restriction map $H^{2}(X;\Z) \to H^{2}(Y;\Z)$ is
        surjective. \label{it:H2Onto}

        \item The manifold $X$ contains a tight surface of positive
        self-intersection number, and a sphere of self-intersection
        $-1$.

        \item The two pieces $X_{1}$ and $X_{2}$ obtained by cutting
        $X$ along $Y$ both have $b^{+}$ positive.
    \end{enumerate}
\end{prop}

\begin{proof}
    By the results of \cite{Eliashberg-Thurston}, the
existence of the foliation implies that
the product manifold
\[
   W = [-1,1]\times Y
\]
carries a symplectic form
$\omega$, weakly compatible with contact structures $\xi_{+}$ and
$\xi_{-}$ on the boundary components $\{1\}\times Y_{0}$ and $\{-1\}
\times Y_{0}$. By Theorem~\ref{thm:Eliash}, we may embed $(W,\omega)$
in a closed symplectic $4$--manifold $(X,\Omega)$.  We can choose $X$ to
satisfy the extra conditions in Lemmas~\ref{lem:noH1},
\ref{lem:H2onto}, \ref{lem:LikeHyp} and \ref{lem:Xtight} above.
This gives the first of the four conditions on $X$. The last condition
is straightforward.
\end{proof}

\subsection{Proof of Theorem~\ref{thm:main}}

Let $Y_{1}$ be the result of $+1$--surgery on a non-trivial knot $K$,
and let $Y_{0}$ be
the manifold with $H_{1}(Y_{0}) = \Z$ obtained by $0$--surgery.
According to Floer's exact triangle \cite{Floer2,Braam-Donaldson}, the
instanton Floer homology group $\HF(Y_{1})$ is isomorphic to the Floer
homology group $\HF(Y_{0})$, where the latter is refers to the group
constructed using the $\SO(3)$ bundle $P \to Y_{0}$ with non-zero
$w_{2}$. We suppose that the knot $K$ contradicts
Theorem~\ref{thm:main}. Then
$\HF(Y_{1})$ is zero, and the exact triangle tells us
that $\HF(Y_{0})$ is zero also.
We therefore have:

\begin{prop}\label{prop:Vanish}
     Suppose $K$ is a counterexample to Theorem~\ref{thm:main}.
    Let $X$ be a smooth closed $4$--manifold containing $Y_{0}$ as a
    separating hypersurface. Suppose that the two pieces $X_{1}$,
    $X_{2}$ obtained by cutting $X$ along $Y_{0}$ both have $b^{+}$
    non-zero. Then the Donaldson polynomial invariant $D_{X}^{w}$ is
    identically zero for any class $w \in H^{2}(X;\Z)$
    whose restriction to $Y_{0}$ is
    non-zero mod $2$.
\end{prop}

\begin{proof}
 When $X$ is decomposed along $Y_{0}$ as in the proposition,
    the value of $D_{X}^{w}(x^{m} h^{n})$ can be expressed as a
    pairing
    \[
               \langle \psi_{X_{1}},\psi_{X_{2}} \rangle,
    \]
    where $\psi_{X_{1}}$ and $\psi_{X_{2}}$ are relative invariants of
    $X_{1}$ and $X_{2}$ taking values in the Fukaya--Floer homology
    group $\HFF(Y_{0},\delta)$ and its dual, where $\delta$ is a
    $1$--cycle in $Y_{0}$ (see \cite{Fukaya,Braam-Donaldson-HFF}). The vanishing of
    $\HF(Y_{0})$ implies the vanishing of $\HFF(Y_{0},\delta)$ also,
    which explains the proposition.
\end{proof}

\subparagraph{\rem} It is possible to avoid the use of the full
exact triangle, and to avoid mentioning any type of Floer homology in
the proof of this proposition. The hypothesis on $K$ means that the
equations for a flat $\SO(3)$ connection on $Y_{0}$ with $w_{2}$
non-zero admit a holonomy-type perturbation (of the sort described in
\cite{Braam-Donaldson}), so that the resulting equations admit no
solutions. (In other language, the Chern--Simons functional has a
holonomy-type perturbation after which it has no critical points.) The
vanishing of the Donaldson invariants for $X$ then follows from a
straightforward degeneration argument.

\medskip
According to \cite{Gabai3}, the manifold $Y_{0}$ has a taut foliation
by oriented 2--dimensional leaves and is not the product manifold
$S^{1}\times S^{2}$ if $K$ is non-trivial. We may apply
Proposition~\ref{prop:Summary} to $Y_{0}$, to embed it in $(X,\Omega)$
satisfying all the conditions in that proposition. 
Being symplectic, the manifold $X$ has $\SW$--simple type and
non-trivial Seiberg--Witten invariants, by the results of
\cite{Taubes1}. The conditions imposed in
Proposition~\ref{prop:Summary} ensure that Corollary~\ref{cor:tightWitten}
applies, so Witten's conjecture, in the form of
Conjecture~\ref{conj:Witten}, holds for $X$. It follows
that the Donaldson invariants $D^{w}_{X}$ are non-trivial, for all
$w$.

However, the $3$--manifold $Y_{0}\subset X$ divides $X$ into two pieces
$X_{1}$ and $X_{2}$, both of which have positive
$b^{+}$. The condition (\ref{it:H2Onto}) of
Proposition~\ref{prop:Summary} allows us to choose a
$w \in H^{2}(X;\Z)$ whose restriction to $Y_{0}$ is the generator.
For this choice of $w$, Proposition~\ref{prop:Vanish} tells us that
$D^{w}_{X}$ is zero. This is a contradiction. \qed

\subsection{Proof of Theorem~\ref{thm:Betti}}

Let $Y$ and $v$ be as in the statement of the theorem.
If the image of the element
$v$ in $\mathrm{Hom}(H_{2}(Y;\Z), \Z/2)$ is zero, then
the result is elementary, for there is an integer lift of $v$ that is
a torsion element of $H^{2}(Y;\Z)$, which implies that there is a flat
$\SO(2)$ bundle on $Y$ with $w_{2}=v$. We therefore
turn to the interesting
case, when $v$ has non-zero pairing with some element of
$H_{2}(Y;\Z)$.

Gabai's theorem
\cite{Gabai} supplies $Y$ with a taut foliation, so we can embed $Y$
as a separating hypersurface in a symplectic $4$--manifold $X$, as in
Proposition~\ref{prop:Summary}. Because the restriction map
on second cohomology is surjective, there is a class
$w\in H^{2}(X;\Z)$ whose restriction to $Y$ becomes $v$ when reduced
mod $2$.

The hypothesis that $v$ has non-zero pairing with some integer class
ensures that there is a
well-defined Floer homology group $\HF^{v}(Y)$ constructed from the
connections with $w_{2}=v$ (see \cite{Donaldson-Book}). If there are
no such
flat connections, then $\HF^{v}(Y)$ is zero, and it follows that
$\bD^{w}_{X}$ is identically zero, as in
Proposition~\ref{prop:Vanish}.
On the other hand,
Conjecture~\ref{conj:Witten} holds for $X$, and we have the same
contradiction as before. \qed


\begin{thebibliography}

\bibitem{Akbulut-McCarthy}
\textbf{Selman Akbulut}, \textbf{John~D McCarthy}, \emph{Casson's invariant for
  oriented homology {$3$}--spheres}, volume~36 of \emph{Mathematical Notes},
  Princeton University Press, Princeton, NJ (1990), an exposition

\bibitem{Bing-Erratum}
\textbf{R\,H Bing}, \emph{Correction to ``{N}ecessary and sufficient conditions
  that a {$3$}--manifold be {$S\sp{3}$}''}, Ann. of Math. 77 (1963) 210

\bibitem{Bing-Martin}
\textbf{R\,H Bing}, \textbf{J\,M Martin}, \emph{Cubes with knotted holes},
  Trans. Amer. Math. Soc. 155 (1971) 217--231

\bibitem{Braam-Donaldson}
\textbf{P\,J Braam}, \textbf{S\,K Donaldson}, \emph{Floer's work on instanton
  homology, knots and surgery}, from: ``The Floer memorial volume'', Progr.
  Math. 133, Birkh\"auser, Basel (1995)  195--256

\bibitem{Braam-Donaldson-HFF}
\textbf{P\,J Braam}, \textbf{S\,K Donaldson}, \emph{Fukaya-{F}loer homology and
  gluing formulae for polynomial invariants}, from: ``The Floer memorial
  volume'', Progr. Math. 133, Birkh\"auser, Basel (1995)  257--281

\bibitem{CGLS}
\textbf{M Culler}, \textbf{C\,McA Gordon}, \textbf{J Luecke}, \textbf{P\,B
  Shalen}, \emph{Dehn surgery on knots}, Ann. of Math. 125 (1987) 237--300

\bibitem{Donaldson-Book}
\textbf{S\,K Donaldson}, \emph{Floer homology groups in {Y}ang--{M}ills theory},
  volume 147 of \emph{Cambridge Tracts in Mathematics}, Cambridge University
  Press, Cambridge (2002), with the assistance of M. Furuta and D. Kotschick

\bibitem{Eliashberg}
\textbf{Y\,M Eliashberg}, \emph{{A few remarks about symplectic filling}},
\href{http://www.maths.warwick.ac.uk/gt/GTVol8/paper6.abs.html}{Geom.\,Topol.\,8\,(2004)} 277-293

\bibitem{Eliashberg-Thurston}
\textbf{Y\,M Eliashberg}, \textbf{W\,P Thurston}, \emph{Confoliations},
  University Lecture Series 13, American Mathematical Society (1998)

\bibitem{Etnyre}
\textbf{J\,B Etnyre}, \emph{{On Symplectic Fillings}}, \agtref4{2004}5{73}{80}

\bibitem{Etnyre-Honda}
\textbf{J\,B Etnyre}, \textbf{K Honda}, \emph{On symplectic cobordisms}, Math.
  Ann. 323 (2002) 31--39

\bibitem{Feehan-Leness}
\textbf{P\,M\,N Feehan}, \textbf{T\,G Leness}, \emph{{A general SO(3)-monopole
  cobordism formula relating Donaldson and Seiberg--Witten invariants}},
  \arxiv{math.DG/0203047}

\bibitem{Floer2}
\textbf{A Floer}, \emph{Instanton homology and {D}ehn surgery}, from: ``The
  Floer memorial volume'', Progr. Math. 133, Birkh\"auser, Basel (1995)  77--97

\bibitem{Fukaya}
\textbf{K Fukaya}, \emph{Floer homology for oriented {$3$}--manifolds}, from:
  ``Aspects of low-dimensional manifolds'', Adv. Stud. Pure Math. 20,
  Kinokuniya, Tokyo (1992)  1--92

\bibitem{Gabai}
\textbf{D Gabai}, \emph{{Foliations and the topology of 3--manifolds}}, J.
  Differential Geom 18 (1983) 445--503

\bibitem{Gabai3}
\textbf{D Gabai}, \emph{Foliations and the topology of $3$--manifolds.
  {I}{I}{I}}, J. Differential Geom. 26 (1987) 479--536

\bibitem{Giroux}
\textbf{E Giroux}, \emph{G\'eom\'etrie de contact: de la dimension trois vers
  les dimensions sup\'erieures}, from: ``Proceedings of the International
  Congress of Mathematicians, Vol. II (Beijing, 2002)'', Higher Ed. Press,
  Beijing (2002)  405--414

\bibitem{Gordon-Luecke}
\textbf{C\,McA Gordon}, \textbf{J Luecke}, \emph{Knots are determined by their
  complements}, J. Amer. Math. Soc. 2 (1989) 371--415

\bibitem{Kirby}
\textbf{R\,C Kirby}, \emph{Problems in low-dimensional topology}, from:
  ``Proceedings of the 1993 Georgia International Topology Conference'',
  (William~H Kazez, editor), AMS/IP Studies in Advanced Mathematics 2, American
  Mathematical Society, Providence, RI (1997)  35--374

\bibitem{KM-Structure}
\textbf{P\,B Kronheimer}, \textbf{T\,S Mrowka}, \emph{Embedded surfaces and the
  structure of {D}onaldson's polynomial invariants}, J. Differential Geom. 41
  (1995) 573--734

\bibitem{Taubes1}
\textbf{C\,H Taubes}, \emph{{The Seiberg--Witten invariants and symplectic
  forms}}, Math. Res. Lett. 1 (1994) 809--822

\bibitem{Witten}
\textbf{E Witten}, \emph{{Monopoles and four-manifolds}}, Math. Res. Lett. 1
  (1994) 769--796

\end{thebibliography}
\end{document}